\documentclass[12pt]{article}
\usepackage{eqsection,subeqnarray,indent,amsfonts,epsf}

\begin{document}

\bibliographystyle{plain}

\newcommand{\be}{\begin{equation}}
\newcommand{\ee}{\end{equation}}
\newcommand{\<}{\langle}
\renewcommand{\>}{\rangle}
\newcommand{\widebar}{\overline}
\def\reff#1{(\protect\ref{#1})}
\def\spose#1{\hbox to 0pt{#1\hss}}
\def\ltapprox{\mathrel{\spose{\lower 3pt\hbox{$\mathchar"218$}}
 \raise 2.0pt\hbox{$\mathchar"13C$}}}
\def\gtapprox{\mathrel{\spose{\lower 3pt\hbox{$\mathchar"218$}}
 \raise 2.0pt\hbox{$\mathchar"13E$}}}
\def\textprime{${}^\prime$}
\def\proof{\par\medskip\noindent{\sc Proof.\ }}
\def\qed{\hbox{\hskip 6pt\vrule width6pt height7pt depth1pt \hskip1pt}\bigskip}
\def\proofof#1{\bigskip\noindent{\sc Proof of #1.\ }}
\def\half{ {1 \over 2} }
\def\third{ {1 \over 3} }
\def\twothird{ {2 \over 3} }
\def\smfrac#1#2{{\textstyle{#1\over #2}}}
\def\smhalf{ \smfrac{1}{2} }
\newcommand{\real}{\mathop{\rm Re}\nolimits}
\renewcommand{\Re}{\mathop{\rm Re}\nolimits}
\newcommand{\imag}{\mathop{\rm Im}\nolimits}
\renewcommand{\Im}{\mathop{\rm Im}\nolimits}
\newcommand{\sgn}{\mathop{\rm sgn}\nolimits}
\newcommand{\tr}{\mathop{\rm tr}\nolimits}
\newcommand{\diag}{\mathop{\rm diag}\nolimits}
\newcommand{\Gal}{\mathop{\rm Gal}\nolimits}
\newcommand{\mycup}{\mathop{\cup}}
\def\hboxscript#1{ {\hbox{\scriptsize\em #1}} }
\def\zhat{ {\widehat{Z}} }
\def\phat{ {\widehat{P}} }
\def\qtilde{ {\widetilde{q}} }
\newcommand{\mod}{\mathop{\rm mod}\nolimits}

\def\Z{{\mathbb Z}}
\def\R{{\mathbb R}}
\def\C{{\mathbb C}}
\def\D{{\mathbb D}}
\def\Q{{\mathbb Q}}

\newtheorem{theorem}{Theorem}[section]
\newtheorem{proposition}[theorem]{Proposition}
\newtheorem{lemma}[theorem]{Lemma}
\newtheorem{corollary}[theorem]{Corollary}
\newtheorem{conjecture}[theorem]{Conjecture}


\newenvironment{sarray}{
          \textfont0=\scriptfont0
          \scriptfont0=\scriptscriptfont0
          \textfont1=\scriptfont1
          \scriptfont1=\scriptscriptfont1
          \textfont2=\scriptfont2
          \scriptfont2=\scriptscriptfont2
          \textfont3=\scriptfont3
          \scriptfont3=\scriptscriptfont3
        \renewcommand{\arraystretch}{0.7}
        \begin{array}{l}}{\end{array}}

\newenvironment{scarray}{
          \textfont0=\scriptfont0
          \scriptfont0=\scriptscriptfont0
          \textfont1=\scriptfont1
          \scriptfont1=\scriptscriptfont1
          \textfont2=\scriptfont2
          \scriptfont2=\scriptscriptfont2
          \textfont3=\scriptfont3
          \scriptfont3=\scriptscriptfont3
        \renewcommand{\arraystretch}{0.7}
        \begin{array}{c}}{\end{array}}


\title{Numerical Computation of
       $\prod\limits_{n=1}^\infty (1 - tx^n)$}
 
\author{
  \\[-1mm]
  {Alan D. Sokal}                   \\
  {\it Department of Physics}       \\
  {\it New York University}         \\
  {\it 4 Washington Place}          \\
  {\it New York, NY 10003 USA}      \\
  {\tt SOKAL@NYU.EDU}          \\
  \\
}
\vspace{0.5cm}
 
\date{December 2, 2002}
\maketitle
\thispagestyle{empty}   

\begin{abstract}
I present and analyze a quadratically convergent algorithm
for computing the infinite product $\prod_{n=1}^\infty (1 - tx^n)$
for arbitrary complex $t$ and $x$ satisfying $|x| < 1$,
based on the identity
$$
   \prod\limits_{n=1}^\infty (1 - tx^n)
   \;=\;  \sum_{m=0}^\infty
          {(-t)^m x^{m(m+1)/2}  \over (1-x)(1-x^2) \cdots (1-x^m)} 
$$
due to Euler.
The efficiency of the algorithm deteriorates as $|x| \uparrow 1$,
but much more slowly than in previous algorithms.
The key lemma is a two-sided bound on the
Dedekind eta function at pure imaginary argument, $\eta(iy)$,
that is sharp at the two endpoints $y=0,\infty$
and is accurate to within 9.1\% over the entire interval $0 < y < \infty$.
\end{abstract}
 
\vspace{0.5cm}
\noindent
{\bf 2000 MATHEMATICS SUBJECT CLASSIFICATION:}
Primary 33F05;
Secondary 05A30, 11F20, 11P82, 33D99, 65D20, 82B23.

\vspace{0.5cm}
 
\clearpage

\section{Introduction}

The function
\be
   R(t,x) \;=\; \prod\limits_{n=1}^\infty (1 - tx^n) \;,
 \label{def_R}
\ee
defined for complex $t$ and $x$ satisfying $|x| < 1$,
was first studied by Euler \cite{Euler}
and has numerous applications in combinatorics, number theory,
analytic-function theory and statistical mechanics.
The case $t=1$ is equivalent to the Dedekind eta function
\be
   \eta(\tau)  \;=\;  e^{\pi i \tau/12} R(1, e^{2\pi i \tau})
   \;,
\ee
which is a modular form \cite{Apostol_90,Knopp_70}
and plays a central role in the
enumeration of partitions \cite{Andrews_98,Apostol_90,Knopp_70}
and sums of squares \cite{Knopp_70}.
The case $t=-1$ is related to $t=1$ via the trivial identity
\be
   R(-1,x)  \;=\;  {R(1,x^2) \over R(1,x)}   \;.
 \label{eq_Rminus1}
\ee
Both of these cases are related to theta functions
\cite{Andrews_99,Bellman_61,Chandrasekharan_85,Rademacher_73}
via the identities
\begin{eqnarray}
   R(1,x) \;\equiv\; \prod_{n=1}^\infty (1-x^n)
       & = & \sum_{m=-\infty}^\infty (-1)^m x^{m(3m+1)/2}
     \label{Euler_pentagonal} \\[2mm]
   R(1,x)^3 \;\equiv\; \prod_{n=1}^\infty (1-x^n)^3
       & = & \sum_{m=0}^\infty (-1)^m (2m+1) x^{m(m+1)/2}
     \\[2mm]
   {R(1,x^2)^2 \over R(1,x)} \;\equiv\;
      \prod_{n=1}^\infty {1-x^{2n} \over 1-x^{2n-1}}
       & = & \sum_{m=0}^\infty x^{m(m+1)/2}
     \\[2mm]
   \hspace*{-1.2cm}
   {R(1,x^2)^5 \over R(1,x)^2 R(1,x^4)^2} \;\equiv\;
      \prod_{n=1}^\infty (1-x^{2n}) (1+x^{2n-1})^2
       & = & \sum_{m=-\infty}^\infty x^{m^2}
     \\[2mm]
   {R(1,x)^2 \over R(1,x^2)} \;\equiv\; \prod_{n=1}^\infty {1-x^n \over 1+x^n}
       & = & \sum_{m=-\infty}^\infty (-1)^m x^{m^2}
\end{eqnarray}
due to Euler, Jacobi and Gauss,
which have spawned a plethora of modern extensions
 \cite{Gordon_61,Macdonald_72,Lepowsky_78,Kac_78,Kac_80,Neher_85,Kohler_90}.
Additional cases of the function $R(t,x)$
arise in the celebrated Rogers--Ramanujan identities
\cite{Andrews_98,Andrews_99}
\begin{eqnarray}
   \prod_{n=0}^\infty {1 \over (1-x^{5n+1})(1-x^{5n+4})}
   & = &  \sum_{m=0}^\infty  {x^{m^2} \over (1-x)(1-x^2) \cdots (1-x^m)}
         \\[2mm]
   \prod_{n=0}^\infty {1 \over (1-x^{5n+2})(1-x^{5n+3})}
   & = &  \sum_{m=0}^\infty  {x^{m(m+1)} \over (1-x)(1-x^2) \cdots (1-x^m)}
\end{eqnarray}
which have numerous combinatorial consequences \cite{Andrews_98}
and which play a key role in Baxter's solution of the hard-hexagon problem
in statistical mechanics \cite{Baxter_81,Baxter_82}.
(See also \cite{Slater_52,Loxton_84,Gasper_90,Andrews_92,Andrews_98}
 for many related identities.)
Finally --- and this was the initial motivation for the current work ---
the cases $t = \pm 1$ and $t = \pm \omega$,
where $\omega$ is a cube root of unity,
arise in Baxter's solution for the chromatic polynomials
of large triangular lattices \cite{Baxter_86,Baxter_87}.
To determine the limiting curves of chromatic roots for these lattices,
it is necessary to compute $R(t,x)$ to high precision for complex $x$,
including points $x$ very near the unit circle \cite{transfer3}.

The numerical computation of $R(t,x)$ clearly becomes delicate
when $|x| \uparrow 1$.
Surprisingly, there seem to be very few treatments of this problem
in the literature
\cite{Slater_54,Slater_66,Gatteschi_69,Allasia_80},
and the algorithms employed there are only linearly convergent;
moreover, these authors (with the exception of Gatteschi \cite{Gatteschi_69})
considered almost exclusively the case of real $t$ and $x$.
My purpose here is to propose and analyze
a quadratically convergent algorithm for computing $R(t,x)$
for arbitrary complex $t$ and $x$ satisfying $|x| < 1$,
based on the identity
\be
   R(t,x)  \;\equiv\;  \prod\limits_{n=1}^\infty (1 - tx^n)
     \;=\;  \sum_{m=0}^\infty
            {(-t)^m x^{m(m+1)/2}  \over (1-x)(1-x^2) \cdots (1-x^m)} 
 \label{eq_qseries}
\ee
due to Euler.\footnote{
   For a proof of \reff{eq_qseries},
   see e.g.\ \cite[p.~19, Corollary 2.2]{Andrews_98},
   \cite[p.~34, Lemma 4(a)]{Knopp_70}
   or \cite[pp.~22--23]{Remmert_98}.
}
In the course of this analysis, I will obtain (Corollary~\ref{cor.R0})
a two-sided bound on $R(1,x)$
for $0 < x < 1$ (and thus on $\eta(iy)$ for $0 < y < \infty$)
that is sharp at the two endpoints $x=0,1$
and is accurate to within 9.1\% over the entire interval;
this bound is perhaps of some modest independent interest.

Of course, for the special case $t=1$
one may employ an even faster algorithm based on
using the modular transformation law for the Dedekind eta function
\cite{Andrews_98,Apostol_90,Chandrasekharan_85,Knopp_70,Rademacher_73}
to move $x$ away from the unit circle,
followed by evaluation of the quadratically convergent sum
\reff{Euler_pentagonal}.\footnote{
   Surprisingly, I have been unable to find in the literature
   any discussion of such an algorithm.
   The details of its implementation --- in particular, how to find
   an appropriate modular transformation ---
   may not be entirely trivial.
}
Moreover, the case $t=-1$ can be reduced to $t=1$ via \reff{eq_Rminus1}.
But for $t \neq \pm 1$ no such identities are known.

The plan of this paper is as follows:
In Section~\ref{sec2} I formulate and prove the properties
of the function $R(t,x)$ that will be needed in the sequel.
In Section~\ref{sec_numerical} I obtain bounds
(both {\em a priori}\/ and {\em a posteriori}\/)
on the rate of convergence of the algorithm defined by \reff{eq_qseries}.
Finally, in Section~\ref{sec4}, I briefly compare this algorithm
to other algorithms that have been proposed
\cite{Slater_54,Slater_66,Gatteschi_69,Allasia_80}
for computing $R(t,x)$.

\clearpage

\section{Properties of \boldmath$R(t,x)$}   \label{sec2}

We shall assume throughout this note that $|x| < 1$,
even if it is not explicitly stated.
Likewise, when we write $x = e^{-\gamma}$,
we shall assume that $\real\gamma > 0$.


\subsection{Elementary properties}

We begin by noting some elementary properties of $R(t,x)$:

\medskip

1) $R(t,x)$ is a jointly analytic function of $t$ and $x$
for $(t,x) \in \C \times \D$, where $\D$ is the open unit disc.
For fixed $x \in \D$, $R(t,x)$ is an entire function of $t$ of order 0,
with simple zeros (for $x \neq 0$) at $t = x^{-1}, x^{-2}, \ldots\;$.

2) By splitting off the first term in the product \reff{def_R}, we obtain
the functional equation
\be
   R(t,x)  \;=\;  (1-tx) \, R(tx,x)   \;,
 \label{eq.functional}
\ee
from which Euler's formula \reff{eq_qseries} can easily be derived
by comparing coefficients of powers of $t$.

3) Let $\omega$ be a primitive $m$th root of unity, and use the identity
$\prod_{j=0}^{m-1} (1- \omega^j z) = 1-z^m$;  we find
\be
   \prod\limits_{j=0}^{m-1}  R(\omega^j t,x)  \;=\;  R(t^m,x^m)
   \;.
\ee
The special case $m=2$, $t=1$ is \reff{eq_Rminus1}.

4) By splitting the product \reff{def_R} according to residue classes
modulo $m$, we obtain
\be
   R(t,x)  \;=\;  \prod\limits_{j=1}^{m} R(tx^{j-m},x^m)
   \;.
 \label{eq_residue_split}
\ee
This formula permits the determination of the asymptotic behavior of $R(t,x)$
as $x$ approaches an $m$th root of unity,
once the asymptotic behavior as $x \to 1$ is known.

5) We have the trivial upper bound
\be
   |R(t,x)|  \;\le\;  R(-T,|x|)
\ee
whenever $|t| \le T$.

6) We have the trivial lower bound
\be
   |R(t,x)|  \;\ge\;  R(T,|x|)
 \label{eq_trivial_lower}
\ee
whenever $|t| \le T \le |x|^{-1}$.
(The condition $T \le |x|^{-1}$ is of course best possible,
 since $R(t,x)$ vanishes at $t = x^{-1}$.)

7) Finally (and most importantly), 
let us take the logarithm (principal branch)
of the defining equation \reff{def_R},
expand $\log(1 - tx^n)$ in Taylor series,
and interchange the absolutely convergent summations;
this yields the useful representation as a Lambert series
\be
   \log R(t,x)  \;=\;
   - \sum_{k=1}^\infty {t^k \over k} \, {x^k \over 1-x^k}
   \;,
 \label{eq_logR}
\ee
valid whenever $|x| < 1$ and $|tx| < 1$.
We shall use this representation repeatedly.

\subsection{Elementary bounds}


Bounding the denominator of \reff{eq_logR}
using $|1 - x^k| \ge 1 - |x|^k \ge 1 - |x|$, we obtain:

\begin{lemma}
   \label{lemma2.1}
Whenever $|x| < 1$ and $|tx| < 1$, we have
\be
   |\log R(t,x)|  \;\le\;  {-\log(1-|tx|)  \over  1-|x|}
\ee
(where the principal branch of the logarithm is taken)
and hence
\be
   (1-|tx|)^{1/(1-|x|)}  \;\le\; |R(t,x)|  \;\le\;  (1-|tx|)^{-1/(1-|x|)}
   \;.
\ee
\end{lemma}

\noindent
This is a crude bound that does {\em not}\/ exhibit the correct
behavior as $|x| \to 1$, but we shall use it as a starting point
for further refinements.

First we need a slight extension of Lemma~\ref{lemma2.1}
for the special case $t=1$.
For $0 \le x < 1$, define
\be
   S(x)  \;=\;  -\log R(1,x)
         \;=\;  \sum_{k=1}^\infty {1 \over k} \, {x^k \over 1-x^k}  \;,
 \label{def_S}
\ee
so that
\begin{eqnarray}
   S'(x)   & = &  \sum_{k=1}^\infty {x^{k-1} \over (1-x^k)^2}
      \label{def_Sprime} \\
   S''(x)  & = &  \sum_{k=1}^\infty
      {(k-1) x^{k-2} + (k+1) x^{2k-2} \over (1-x^k)^3}
      \label{def_Sdoubleprime}
\end{eqnarray}
By using $1 \ge 1-x^k \ge 1-x$ in the denominator,
we obtain the trivial bounds:

\begin{lemma}
   \label{lemma2.2}
For $0 \le x < 1$, we have
\begin{eqnarray}
   -\log(1-x)  \;\le\;  S(x)  & \le &
       {-\log(1-x) \over 1-x}  \;\le\;  {x \over (1-x)^2}
      \label{bound_S}   \\[2mm]
   1  \;\le\; {1 \over 1-x}  \;\le\;  S'(x)  & \le &
       {1 \over (1-x)^3}
      \label{bound_Sprime}   \\[2mm]
   3  \;\le\; {1 \over (1-x)^2} + {2-x^2 \over (1-x^2)^2} \;\le\;  S''(x)
           & \le &
       {1 \over (1-x)^5} + {2-x^2 \over (1-x)^3 (1-x^2)^2}  \;\le\;
       {3 \over (1-x)^5}
      \nonumber \\ \label{bound_Sdoubleprime}
\end{eqnarray}
\end{lemma}

\subsection{Case \boldmath$t=1$}

Now we improve these bounds by using a deep fact:
the transformation properties of the Dedekind eta function
under the modular group
\cite{Andrews_98,Apostol_90,Chandrasekharan_85,Knopp_70,Rademacher_73}.
All we need, in fact, is a special case of the modular transformation law,
namely the one for inversion $\tau \to -1/\tau$:
\be
   R(1, e^{-2\pi z})  \;=\;
   z^{-1/2} \exp\!\left( {\pi z \over 12} - {\pi \over 12z} \right)
   R(1, e^{-2\pi/z})
 \label{eq_modular}
\ee
for $\real z > 0$.\footnote{
   There are a number of proofs of \reff{eq_modular}.
   The simplest uses the Poisson summation formula
   applied to Euler's pentagonal number theorem \reff{Euler_pentagonal}
   \cite[Section 3.3]{Knopp_70}.
   Another proof, due to Siegel, uses the Cauchy integral formula
   \cite[Section 3.2]{Apostol_90} \cite[Section VIII.3]{Chandrasekharan_85}.
   Proofs of the full modular transformation law
   are given in \cite[Sections 3.3--3.6 and pp.~190--195]{Apostol_90},
   \cite[Sections 3.1--3.3 and 4.1--4.2]{Knopp_70},
   \cite[Chapter 9]{Rademacher_73},
   and \cite[pp.~82--85]{Andrews_98}.
}
This allows us to control the behavior near $x=1$ ($z \to 0$)
in terms of the (trivial) behavior near $x=0$ ($z \to +\infty$).\footnote{
   Using the full modular transformation law, one can control
   in an analogous way the behavior of $R(1,x)$
   near any point $x = e^{2\pi i h/k}$ ($h,k \in \Z$)
   of the unit circle:  see e.g.\ \cite[Chapter 5]{Andrews_98}.
}
Indeed, from \reff{eq_modular} and the regularity of $R(1,x)$ near $x=0$,
one immediately deduces the sharp asymptotic formula
\be
   \log R(1,e^{-\gamma})  \;=\;
   -\, {\pi^2 \over 6\gamma} \,-\, {1 \over 2} \log\gamma \,+\,
        {1 \over 2} \log(2\pi) \,+\, {\gamma \over 24}  \,+\,
        O(e^{-4\pi^2/\gamma})
 \label{eq.sharp.R1}
\ee
as $\gamma \to 0$;
moreover, an explicit quantitative bound on the $O(e^{-4\pi^2/\gamma})$ term
can easily be extracted from Lemma~\ref{lemma2.1}.

For later applications we need also a quantitative error bound
valid for real $\gamma$ in the entire interval $0 < \gamma < \infty$.
Let us define
\be
   f(z)  \;=\;  \log R(1, e^{-2\pi z})  \,+\,  {\pi \over 12z}
                   \,-\, {1 \over 4} \log\!\left(1 + {1 \over z^2} \right)
   \;,
 \label{def_f}
\ee
so that
\be
  f(1/z)  \;=\;  \log R(1, e^{-2\pi/z})  \,+\,  {\pi z \over 12}
                   \,-\, {1 \over 4} \log(1 +  z^2)
  \;.
 \label{eq_f_1overz}
\ee
Then the transformation law \reff{eq_modular} tells us immediately that
$f(z) = f(1/z)$, and indeed we have:

\begin{proposition}
   \label{prop_f}
For $0 < z < \infty$, we have:
\begin{itemize}
   \item[(a)]  $f(z) = f(1/z)$
   \item[(b)]  $\lim\limits_{z \downarrow 0} f(z) = 0$
       and $\lim\limits_{z \to +\infty} f(z) = 0$
   \item[(c)]  $0 < f(z) \le f(1) =
                   {\pi \over 6} - {1 \over 4} \log 2 + \log \eta(i)
                   \approx 0.0866399$
   \item[(d)]  $f'(0) = \pi/12$, $f'(z) > 0$ for $0 < z < 1$,
               $f'(1) = 0$, and $f'(z) < 0$ for $z > 1$
   \item[(e)]  $f''(0) = -1/2$, and there exists $z_* > 1$
       such that $f''(z) < 0$ for $0 < z < z_*$ and $f''(z_*) = 0$.
\end{itemize}
\end{proposition}    

\proof
We have already proven that $f(z) = f(1/z)$,
so we can use \reff{def_f} and \reff{eq_f_1overz}
interchangeably as formulae for $f(z)$.
The limiting values of $f$ and its derivatives at $z=0$
can be read off \reff{eq_f_1overz}.

To prove $f(z) > 0$ for $0 < z < \infty$,
it suffices to prove it for $0 < z \le 1$.
Using \reff{eq_f_1overz}, we make the following crude bounds:
\begin{subeqnarray}
   f(z)  & = &  -S(e^{-2\pi/z}) \,+\, {\pi z \over 12}
                   \,-\, {1 \over 4} \log(1 +  z^2)        \\[2mm]
         & \ge &  - {e^{-2\pi/z} \over (1-e^{-2\pi/z})^2}
            \,+\, {\pi z \over 12} \,-\, {1 \over 4} \log(1 +  z^2)   \\[2mm]
         & \ge &  - {e^{-2\pi/z} \over (1-e^{-2\pi/z})^2}
            \,+\, \left( {\pi \over 12} - {1 \over 4} \right) z
\end{subeqnarray}
where we have used \reff{bound_S} and the fact that $0 < z \le 1$.
So we need only show that
\be
   {x \over (1-x)^2}  \;<\;
       \left( {\pi \over 12} - {1 \over 4} \right)
       \left( - {2\pi \over \log x} \right)
\ee
for $0 < x \le e^{-2\pi}$.
But $-x \log x/(1-x)^2$ is an increasing function of $x$
for $0 < x \le 1$,
and its value at $x = e^{-2\pi}$ is
$2 \pi e^{-2\pi}/(1-e^{-2\pi})^2 \approx 0.011777 <
 (\pi/12 - 1/4)(2\pi) \approx 0.074138$.

Next let us prove that there exists $\epsilon > 0$
such that $f''(z) < 0$ for $0 < z < 1+\epsilon$.
Differentiating \reff{eq_f_1overz} twice with respect to $z$, we obtain
\be
   f''(z)  \;=\;
      \left(\! -{4\pi^2 \over z^4} + {4\pi \over z^3} \right)
         e^{-2\pi/z} S'(e^{-2\pi/z})
      \,-\,
      {4\pi^2 \over z^4} e^{-4\pi/z} S''(e^{-2\pi/z})
      \,-\,
      {1-z^2 \over 2(1+z^2)^2}
   \;.
 \label{eq_fdoubleprime}
\ee
From \reff{bound_Sprime}/\reff{bound_Sdoubleprime}
we have $S'(x) \ge 1$ and $S''(x) \ge 3$,
so the first two terms in \reff{eq_fdoubleprime} are $<0$ for $0<z<\pi$,
and the third term is $\le 0$ for $0 < z \le 1$.
This proves the claim.

We have just proven that $f'(z)$ is a strictly decreasing function of $z$
on $0 < z < 1+\epsilon$.
{}From $f(z)=f(1/z)$ it follows that $f'(1) = 0$.
Therefore $f'(z) > 0$ for $0 < z < 1$;
by $f(z)=f(1/z)$ it follows that $f'(z) < 0$ for $z > 1$;
and thus $f(z) \le f(1)$ for all $z$.

Finally, it is not possible that $f''(z) < 0$ for all $z$,
as this would imply that $f(z) < 0$ for some $z \in (1,\infty)$.
So we can define $z_* > 1$ to be the smallest $z$ such that $f''(z)=0$.
\qed

\begin{figure}[t]
\begin{center}
\epsfxsize=0.45\textwidth
\leavevmode\epsffile{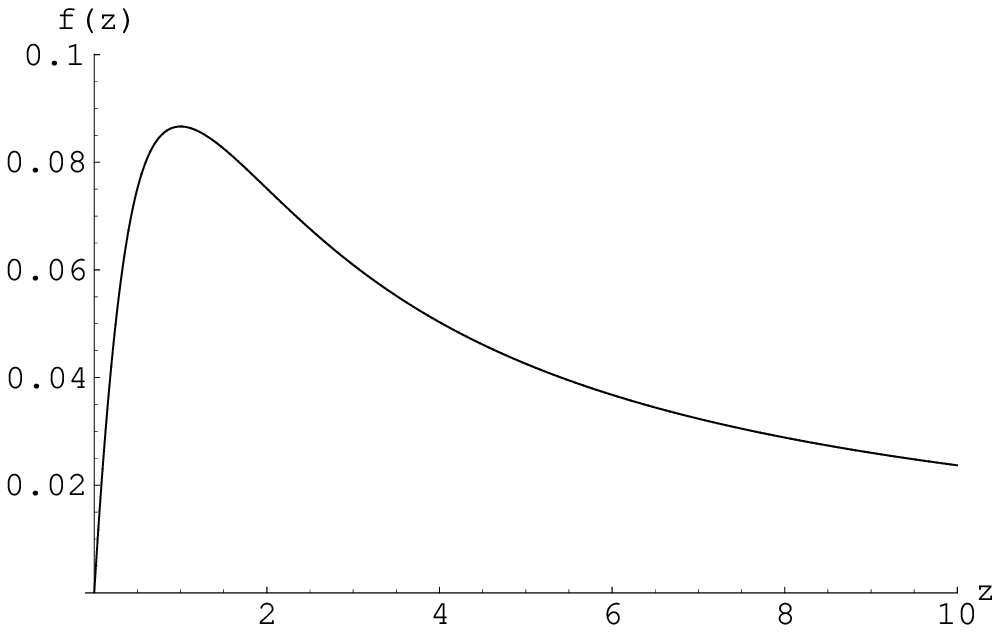}
\hspace{1cm}
\epsfxsize=0.45\textwidth
\epsffile{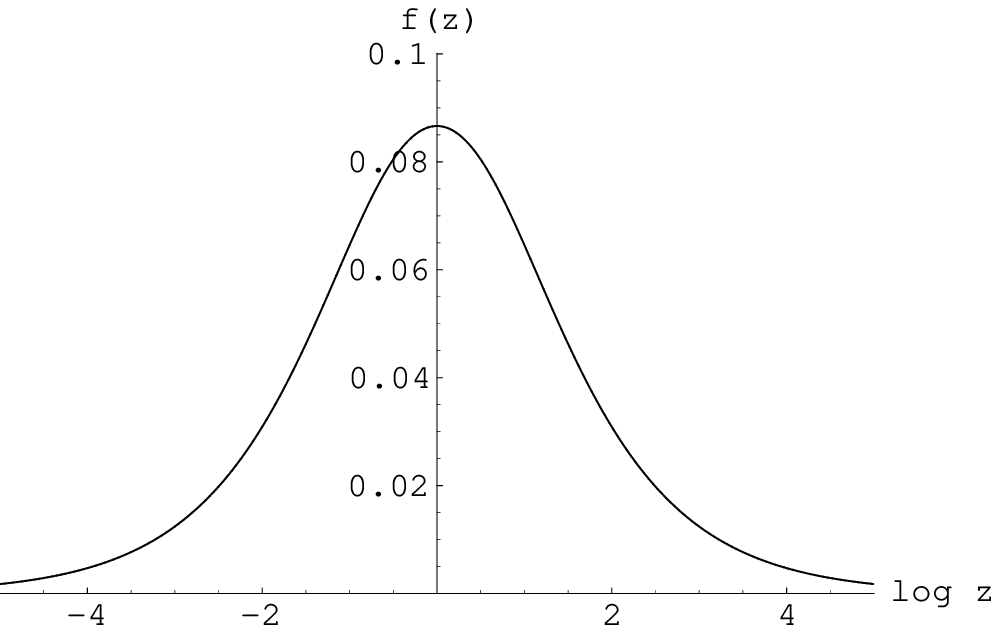} \\
\end{center}
\caption{
   Graphs of $f(z)$ versus $z$ and $\log z$.
}
  \label{fig_f}
\end{figure}

\noindent
{\bf Remark.}
Numerical calculations show that $f''$ has a {\em unique}\/ zero,
which is located at $z_* \approx 1.974174$.
But we shall not bother to prove this.
Graphs of $f(z)$ versus $z$ and $\log z$ are shown in Figure~\ref{fig_f};
the latter shows the $z \leftrightarrow 1/z$ symmetry more clearly.

\bigskip

Proposition~\ref{prop_f} can be rephrased by defining
\be
   R_0(1,x)  \;=\;
      e^{\pi^2/(6\log x)}
      \left( 1 + {4\pi^2 \over (\log x)^2} \right)^{\! 1/4}
   \;,
\ee
which we interpret as an ``approximate'' version of $R(1,x)$.
We then have:

\begin{corollary}
 \label{cor.R0}
For $0 < x < 1$,
\be
   e^{\pi^2/(6\log x)}  \;<\;
   R_0(1,x)  \;<\;  R(1,x)  \;\le\;  C R_0(1,x)
\ee
where $C = e^{f(1)} = e^{\pi/6} 2^{-1/4} \eta(i) \approx 1.090504$.
\end{corollary}

\noindent
In other words, we have a two-sided bound on $R(1,x)$,
in which the lower bound is sharp at the two endpoints $x=0,1$
and is accurate to within 9.1\% over the entire interval $0 < x < 1$.
We shall frequently use the lower bound of Corollary~\ref{cor.R0}
in the form
\be
   R(1,e^{-\gamma})  \;\ge\;
   e^{-\pi^2/6\gamma} \left( 1 + {4\pi^2 \over \gamma^2} \right)^{\! 1/4}
   \;\ge\;  e^{-\pi^2/6\gamma}
 \label{eq.cor.R0.gamma}
\ee
for $\gamma > 0$.

\subsection{Case \boldmath$t=-1$}

We can now handle the case $t=-1$
by using \reff{eq_Rminus1} to relate it to $t=1$.
{}From \reff{eq.sharp.R1} and \reff{eq_Rminus1}
we obtain the sharp asymptotic formula
\be
   \log R(-1,e^{-\gamma})  \;=\;
   {\pi^2 \over 12\gamma} \,-\, {1 \over 2} \log 2
        \,+\, {\gamma \over 24}  \,+\, O(e^{-\pi^2/\gamma})
 \label{eq.sharp.Rminus1}
\ee
as $\gamma \to 0$,
where again a quantitative bound on the $O(e^{-\pi^2/\gamma})$ term
can easily be extracted from Lemma~\ref{lemma2.1}. 
Moreover, we can obtain a quantitative error bound
valid for real $z$ in the entire interval $0 < z < \infty$.
Let us define
\begin{subeqnarray}
   g(z)  & = &  f(\sqrt{2} z)  \,-\,  f(z/\sqrt{2})     \\[2mm]
         & = &  \log R(-1, e^{-\sqrt{2} \pi z})  \,-\,
                  {\pi \over 12 \sqrt{2} z}  \,+\,
                  {1 \over 4} \log\!\left( {z^2 + 2 \over z^2 + \half} \right)
   \;.
 \label{def_g}
\end{subeqnarray}
It follows immediately from Proposition~\ref{prop_f} that:

\begin{proposition}
   \label{prop_g}
For $0 < z < \infty$, we have:
\begin{itemize}
   \item[(a)]  $g(z) = -g(1/z)$
   \item[(b)]  $\lim\limits_{z \downarrow 0} g(z) = 0$
       and $\lim\limits_{z \to +\infty} g(z) = 0$
   \item[(c)]  $g'(z) < 0$ for $1/\sqrt{2} \le z \le \sqrt{2}$
   \item[(d)]  $g(z) > 0$ for $1/\sqrt{2} \le z < 1$,
               $g(1) = 0$, and
               $g(z) > 0$ for $1 < z \le \sqrt{2}$
   \item[(e)]  $|g(z)| \le f(1) \approx 0.0866399$ for $0 < z < \infty$
\end{itemize}
\end{proposition}        

\begin{figure}[t]
\begin{center}
\epsfxsize=0.45\textwidth
\leavevmode\epsffile{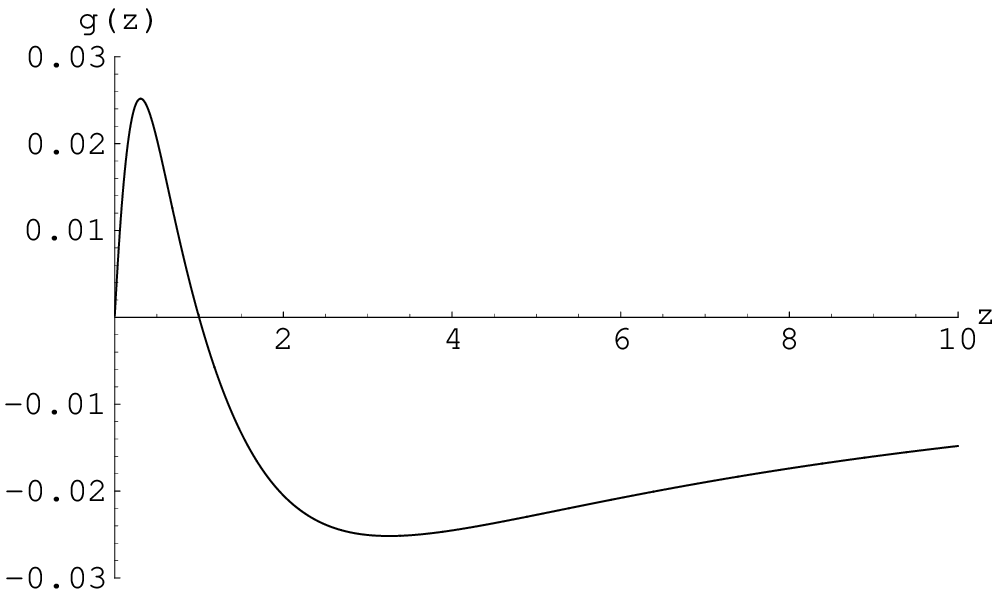}
\hspace{1cm}
\epsfxsize=0.45\textwidth
\epsffile{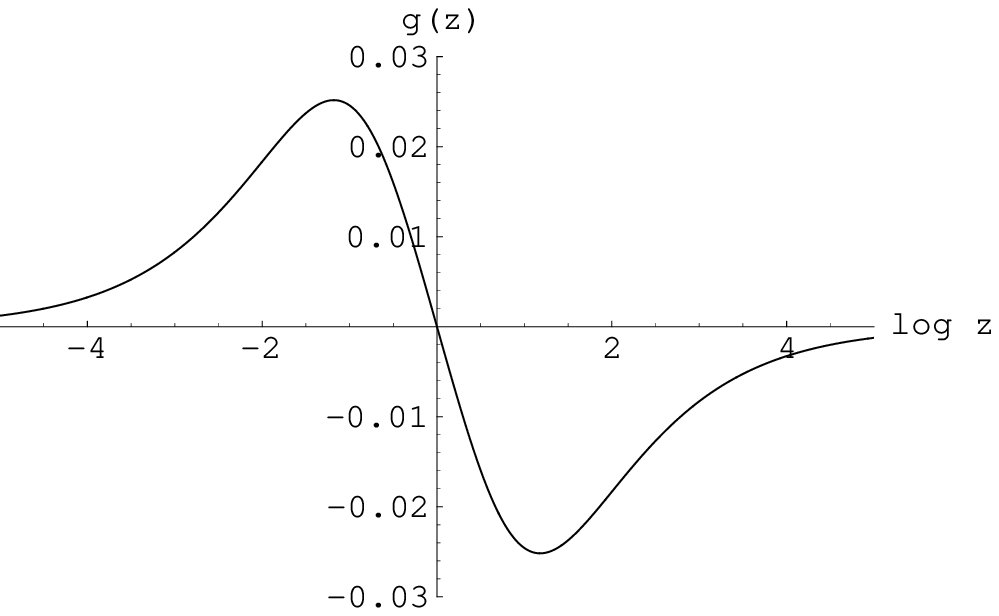} \\
\end{center}
\caption{
   Graphs of $g(z)$ versus $z$ and $\log z$.
}
  \label{fig_g}
\end{figure}

\noindent
{\bf Remark.}
Numerical calculations show that $g'$ vanishes when (and only when)
$\pm \log z \approx 1.180158$,
i.e.\ $z$ or $1/z \approx 3.254889$,
and that the maximum value of $|g(z)|$ is $\approx 0.0251707$.
It follows that $R(-1,x)$ differs from
\be
   R_0(-1,x)  \;\equiv\;
   {R_0(1,x^2) \over R_0(1,x)}
   \;=\;
   e^{-\pi^2/(12\log x)}
      \left( {  1 + {\pi^2 \over (\log x)^2}
                \over
                1 + {4\pi^2 \over (\log x)^2}  } \right)^{\! 1/4}
\ee
by less than 2.6\% over the entire interval $0 < x < 1$.
Graphs of $g(z)$ versus $z$ and $\log z$ are shown in Figure~\ref{fig_g}.

\subsection{Asymptotics of \boldmath$R(t,x)$ for General \boldmath$t$}

Finally, let us discuss briefly the asymptotics of $R(t,x)$ as $x \to 1$
when $t$ is fixed with $|t| < 1$
(or more generally varies within a compact subset of the open unit disc).
Let us write $x = e^{-\gamma}$ with $\real\gamma > 0$
and study the behavior as $\gamma \to 0$,
using the representation \reff{eq_logR}.
We have
\be
   {x^k \over 1-x^k}  \;=\;
   {1 \over e^{k\gamma} - 1}  \;=\;
   \sum\limits_{m=0}^\infty {B_m \over m!} (k\gamma)^{m-1}
 \label{eq_bernoulli}
\ee
where $B_m$ is the $m$th Bernoulli number;
this series is absolutely convergent for $|\gamma| < 2\pi/k$.\footnote{
   See e.g.\ \cite[equation (6.81)]{Graham_94}.
}
Inserting this into \reff{eq_logR}
and {\em formally}\/ interchanging the order of summation,
we obtain\footnote{
   See also \cite[p.~58, exercise 2]{deBruijn_81}
   and \cite[Theorem 4]{McIntosh_99} for this formula.
}
\be
   -\log R(t,e^{-\gamma})  \;\hbox{``=''}\;
   \sum\limits_{m=0}^\infty {B_m \over m!} \, {\rm Li}_{2-m}(t) \, \gamma^{m-1}
 \label{eq_seriesRt}
\ee
where
\be
   {\rm Li}_p(t)  \;\equiv\;  \sum_{k=1}^\infty {t^k \over k^p}
 \label{def_fp}
\ee
is the polylogarithm function \cite{Lewin_81}.
However, because the radius of convergence of \reff{eq_bernoulli}
is nonuniform in $k$ and tends to zero as $k \to \infty$,
it is reasonable to expect that the series \reff{eq_seriesRt}
is {\em not}\/ convergent but is only asymptotic.
One further expects that this asymptotic expansion should hold uniformly
as $t$ varies within a compact subset of the open unit disc.
All these expectations are true \cite{Sokal_dedekind}.
What is perhaps more surprising is that the expansion \reff{eq_seriesRt}
holds also for $t$ {\em on}\/ the unit circle, except at the point $t=1$.
Indeed, under suitable restrictions on $\arg\gamma$
it holds in a much larger domain of the complex $t$-plane,
which in the most favorable case ($\gamma$ real and positive)
encompasses the entire complex $t$-plane
except for a cut along $[1,\infty)$.
These results will be reported elsewhere \cite{Sokal_dedekind}.
For real $\gamma > 0$ and $0 < t < 1$,
the expansion \reff{eq_seriesRt} was proven
some years ago by Moak \cite[Theorem 3]{Moak_84}.\footnote{
   Equation (4.2) of \cite{Moak_84} contains a misprint:
   there should be a minus sign before the integral.
   Correspondingly, in equation (4.3), the minus sign before the integral
   should be a plus sign.
}
For real $\gamma > 0$ and $t \le 1$,
the expansion \reff{eq_seriesRt} and some generalizations thereof
have recently been proven by McIntosh \cite{McIntosh_99}.
For real $\gamma > 0$ and $t \in \C \setminus [1,\infty)$,
the expansion \reff{eq_seriesRt} has been proven by
Prellberg \cite[Lemma 3.2]{Prellberg_95}.
All these works use the Euler--Maclaurin sum formula.
Our approach \cite{Sokal_dedekind}, by contrast, uses complex integration.

For real $\gamma > 0$ and $0 \le t \le 1$,
we can use the method just sketched
to obtain a two-sided bound on $R(t,e^{-\gamma})$
that incorporates the first two terms of the expansion \reff{eq_seriesRt}.
For $z > 0$ we have the elementary inequalities\footnote{
   The first two inequalities can be derived from $\tanh(z/2) \le z/2$;
   the third can be derived from $z/2 \le \sinh(z/2)$;
   and the fourth is trivial.
   Note that all of these bounds, except the last,
   capture the first two terms of the Laurent series for $1/(e^z - 1)$
   around $z=0$.
}
\be
   {1 \over z} \,-\, {1 \over 2}
   \;\le\;
   e^{-z} \left( {1 \over z} \,+\, {1 \over 2} \right)
   \;\le\;
   {1 \over e^z -1}
   \;\le\;
   {e^{-z/2} \over z}
   \;\le\;
   {1 \over z}
   \;.
\ee
Setting $z=k\gamma$ and inserting these bounds into \reff{eq_logR},
we obtain:

\begin{proposition}
 \label{prop_gamma_real}
For $0 \le t \le 1$ and $\gamma > 0$, we have
\begin{subeqnarray}
 \hspace*{-30mm}  
   -\log R(t,e^{-\gamma})
        & \le &
        \gamma^{-1} {\rm Li}_2(t e^{-\gamma/2})
 \slabel{eq_upper_minusRt_a} \\[1mm]
        & \le &
        \gamma^{-1} {\rm Li}_2(t)
\end{subeqnarray}
and
\begin{subeqnarray}
 \hspace*{4mm}  
   -\log R(t,e^{-\gamma})   & \ge &
         \gamma^{-1} {\rm Li}_2(t e^{-\gamma})
               \,-\, \smhalf \log(1 - t e^{-\gamma})
 \slabel{eq_lower_minusRt_a} \\[1mm]
     & \ge &
         \gamma^{-1} {\rm Li}_2(t)  \,+\, \smhalf \log(1-t)
 \slabel{eq_lower_minusRt_b}
\end{subeqnarray}
\end{proposition}

\noindent
It is worth remarking that, even for $t=1$,
the bounds \reff{eq_upper_minusRt_a} and \reff{eq_lower_minusRt_a}
capture the first two terms of the asymptotic expansion \reff{eq.sharp.R1},
i.e.\ they get the correct $\log\gamma$ term.

One application of Proposition~\ref{prop_gamma_real}
is to bounding the partial product
\be
   \prod\limits_{k=1}^n (1-x^k)  \;=\;  {R(1,x) \over R(x^n,x)}
\ee
when $0 < x=e^{-\gamma} < 1$
(and we will usually take $n$ to be of order $1/\gamma$).
Inserting the lower bound \reff{eq.cor.R0.gamma} on $R(1,x)$
and the upper bound \reff{eq_lower_minusRt_b} on $R(x^n,x)$,
we obtain:

\begin{corollary}
   \label{cor_gamma_real}
Let $\gamma > 0$. Then
\be
   \prod\limits_{k=1}^n (1-e^{-k\gamma})  \;\ge\;
   \exp\!\left[ {{\rm Li}_2(e^{-n\gamma}) - \pi^2/6  \over  \gamma}  \right] \,
      (1-e^{-n\gamma})^{1/2}   
   \;.
 \label{eq.cor_gamma_real_1}
\ee
In particular, for $n \le (\log 2)/\gamma$ we have
\be
   \prod\limits_{k=1}^n (1-e^{-k\gamma})  \;\ge\;
   \exp\!\left[ {-(\log 2)^2/2 - \pi^2/12  \over  \gamma}  \right] \,
      (1-e^{-n\gamma})^{1/2}
   \;.
 \label{eq.cor_gamma_real_2}
\ee
\end{corollary}

\noindent
Here \reff{eq.cor_gamma_real_2} follows from \reff{eq.cor_gamma_real_1}
and the well-known fact \cite{Lewin_81,Loxton_84}
\be
   {\rm Li}_2(1/2)  \;=\;  {\pi^2 \over 12} \,-\, {(\log 2)^2 \over 2}  \;.
 \label{eq.dilog.half}
\ee
McIntosh \cite{McIntosh_99} has recently obtained
a complete asymptotic expansion of the partial product
$\prod_{k=1}^n (1- te^{-k\gamma})$ for $n = \mu/\gamma$
($\mu$ fixed, real $\gamma \downarrow 0$)
and either $t=1$ or $t < 1$.

\section{Numerical Computation of \boldmath$R(t,x)$}  \label{sec_numerical}

In this section we discuss the use of Euler's formula
\be
   R(t,x)  \;=\;  \sum_{n=0}^\infty
            {(-t)^n x^{n(n+1)/2}  \over (1-x)(1-x^2) \cdots (1-x^n)}
 \label{eq_qseries_bis}
\ee
to compute $R(t,x)$ for complex $t$ and $x$ satisfying $|x| < 1$.
We shall give two types of bounds on the error committed by
truncating the series \reff{eq_qseries_bis}:
\begin{itemize}
   \item[(a)] an {\em a priori}\/ bound in terms of $|t|$ and $|x|$ alone; and
   \item[(b)] an {\em a posteriori}\/ bound, based on the actual computed
       value of the last included term.
\end{itemize}
We shall also give some guidance about the needed numerical precision
in intermediate stages of the calculation, by comparing the largest
term in the sum to the final answer.

We use the following definitions:
\begin{itemize}
   \item The $n$th term:
     $a_n  \,=\,  {\displaystyle
                   {(-t)^n x^{n(n+1)/2}  \over (1-x)(1-x^2) \cdots (1-x^n)}
                  }$
   \item The partial sum after $N-1$ terms:
     $S_N  \,=\,  \sum\limits_{n=0}^{N-1}  a_n$
   \item The remainder after $N-1$ terms:
     $R_N  \,=\,  \sum\limits_{n=N}^{\infty}  a_n$
   \item The absolute error after $N-1$ terms:
     $\Delta_N = |R_N|$
   \item The relative error after $N-1$ terms:
     $\delta_N = |R_N/R(t,x)|$
   \item The modified relative error after $N-1$ terms:
     $\delta'_N = |R_N/S_N|$
\end{itemize}
Clearly $\delta'_N/(1+\delta'_N) \le \delta_N \le \delta'_N/(1-\delta'_N)$,
so the two types of relative error are essentially indistinguishable when
$\delta_N, \delta'_N \ll 1$.

\begin{lemma}
   \label{lemma3.1}
If $|x| < 1$ and $|t| \, |x|^{N+1} < 1$, then
\be
   \sum\limits_{n=N}^\infty |t^n x^{n(n+1)/2}|
   \;\le\;
   {|t|^N \, |x|^{N(N+1)/2} \over 1 \,-\, |t| \, |x|^{N+1}}
   \;.
\ee
\end{lemma}

\proof
Bound the sum by a geometric series, using
\be
   \left| {t^{n+1} x^{(n+1)(n+2)/2} \over t^n x^{n(n+1)/2} } \right|
   \;=\; |t| \, |x^{n+1}|  \;\le\;  |t| \, |x|^{N+1}
\ee
for $n \ge N$.
\qed

\begin{lemma}
   \label{lemma3.2}
If $|x| \le e^{-\gamma}$ ($\gamma > 0$), then
\be
   \left| \prod\limits_{k=1}^n (1-x^k) \right|
   \;\ge\;
   \prod\limits_{k=1}^n (1-|x|^k)
   \;\ge\;
   \prod\limits_{k=1}^\infty (1-|x|^k)
   \;\equiv\;
   R(1,|x|)
   \;\ge\;
   e^{-\pi^2/6\gamma}
   \;.
\ee
\end{lemma}

\proof
An immediate consequence of Corollary~\ref{cor.R0}.
\qed

\noindent
{\bf Remark.}  An improved bound on the partial product
$\prod_{k=1}^n (1-x^k)$ can be obtained
from Corollary~\ref{cor_gamma_real};
it is advantageous when $n\gamma \not\gg 1$.

\begin{proposition}
   \label{prop3.3}
Suppose that $|x| \le e^{-\gamma}$ with $\gamma > 0$.
\begin{itemize}
   \item[(a)]  If $|t| < e^{(N+1)\gamma}$, then
        $\displaystyle \Delta_N  \,\equiv\,
        \left| \sum\limits_{n=N}^{\infty}  a_n \right|
        \,\le\,
        {|t|^N \, e^{\pi^2/6\gamma - N(N+1)\gamma/2}
         \over
         1 - |t| e^{-(N+1)\gamma}
        }\,$.
   \item[(b)]  If $|t| \le 1$, then
        $\displaystyle \delta_N  \,\equiv\,
        {\left| \sum\limits_{n=N}^{\infty}  a_n \right|  \over |R(t,x)|}
        \,\le\,
        {e^{\pi^2/3\gamma - N(N+1)\gamma/2}  \over
         1 - e^{-(N+1)\gamma}
        }\,$.
   \item[(c)]  If $|t| < e^\gamma$, then
        $\displaystyle \delta_N  \,\equiv\,
        {\left| \sum\limits_{n=N}^{\infty}  a_n \right|  \over |R(t,x)|}
        \,\le\,
        {e^{\pi^2/3\gamma - N(N+1)\gamma/2}  \over
         1 - e^{-(N+1)\gamma}
        }
        \; {1 \over 1 - |t| e^{-\gamma}}\,$.
\end{itemize}
\end{proposition}

\proof
(a) is an immediate consequence of Lemmas~\ref{lemma3.1} and \ref{lemma3.2}.
(b) follows from (a) together with the bound
$|R(t,x)| \ge R(1,|x|) \ge e^{-\pi^2/6\gamma}$ from
\reff{eq_trivial_lower} and Corollary~\ref{cor.R0}.
(c) follows from (b) and \reff{eq.functional}.
\qed

\begin{corollary}
   \label{cor3.4}
Let $K \ge 0$,
and suppose that $|t| \le 1$ and $|x| \le e^{-\gamma}$ ($\gamma > 0$).
\begin{itemize}
   \item[(a)]  If ${\displaystyle N \,\ge\,
                  \sqrt{ {\pi^2 \over 3\gamma^2} \,+\, {2K \over \gamma} } }$,
        then $\Delta_N \le e^{-K}$.
   \item[(b)]  If ${\displaystyle N \,\ge\,
                  \sqrt{ {2\pi^2 \over 3\gamma^2} \,+\, {2K \over \gamma} } }$,
        then $\delta_N \le e^{-K}$.
\end{itemize}
\end{corollary}

\proof
Since $K \ge 0$, we have $N\gamma \ge \pi/\sqrt{3}$ and hence
\be
   {e^{-N\gamma/2} \over 1 - e^{-(N+1)\gamma}}  \;\le\;
   {e^{-\pi/2\sqrt{3}}  \over  1 - e^{-\pi/\sqrt{3}}}  \;\approx\;
   0.482426  \;<\;  1  \;.
\ee
Now $\pi^2/6\gamma - N^2 \gamma/2 \le -K$ in case (a),
and $\pi^2/3\gamma - N^2 \gamma/2 \le -K$ in case (b).
The result then follows from Proposition~\ref{prop3.3}(a,b).
\qed

Please note that the bound in Proposition~\ref{prop3.3}(a)
is asymptotically within 9.1\% of being sharp when $0 < x=e^{-\gamma} < 1$
and $N \gg 1/\gamma$
(and in this case is moreover asymptotically sharp as $\gamma \downarrow 0$);
but it is overly pessimistic in other cases,
because the denominator \mbox{$(1-x)(1-x^2)\cdots(1-x^n)$}
is not really as small as Lemma~\ref{lemma3.2} says it could be.
Likewise, the bound in  Proposition~\ref{prop3.3}(b)
is asymptotically (almost-)sharp when, in addition to the above conditions,
we have $t=1$;
but it is overly pessimistic in other cases,
because $|R(t,x)|$ is not really as small as
the bound $|R(t,x)| \ge R(1,|x|)$ says it could be.

It is thus of some value to provide an {\em a posteriori}\/ bound
on the truncation error that is more realistic, when $x \notin (0,1)$,
than the {\em a priori}\/ bound;
such a bound can be used a stopping criterion in the numerical algorithm.
We need the following elementary observation:

\begin{lemma}
   \label{lemma3.5}
If $|x| \le e^{-\gamma}$  with $\gamma > 0$, then
\be
   \left| {a_n \over a_{n-1}} \right|
   \;=\;
   {|t| \, |x|^n  \over  |1-x^n|}
   \;\le\;
   {|t| \, |x|^n  \over  1-|x|^n}
   \;\le\;
   {|t| \, e^{-n\gamma}  \over  1 - e^{-n\gamma}}
   \;.
 \label{eq.lemma3.5}
\ee
\end{lemma}   

\noindent
Lemma~\ref{lemma3.5} tells us that, at least for $0 < x < 1$,
the terms $a_n$ {\em increase}\/ in magnitude until $|x|^n \approx 1/(1+|t|)$,
i.e.\ $n \approx [\log(1+|t|)]/\gamma$, and then decrease.
(For general complex $x$, the terms will {\em sometimes}\/ increase
 up to this point, i.e.\ for those $n$ for which
 $\arg x^n \approx 0$ mod $2\pi$.  How often this occurs depends on
 the Diophantine properties of $\arg x$.)
We can use Lemma~\ref{lemma3.5} to bound the tail of the sum
by a geometric series:

\begin{proposition}
   \label{prop3.6}
Suppose that $|x| \le e^{-\gamma}$ ($\gamma > 0$)
and $N > [\log (1+|t|)]/\gamma$.  Then:
\begin{itemize}
   \item[(a)]  $\displaystyle \Delta_N  \,\equiv\,
        \left| \sum\limits_{n=N}^{\infty}  a_n \right|
        \,\le\,
        |a_{N-1}| \, {|t| e^{-N\gamma} \over 1 - (1+|t|)e^{-N\gamma}}$
   \item[(b)]  $\displaystyle \delta'_N  \,\equiv\,
        {\left| \sum\limits_{n=N}^{\infty}  a_n \right|  \over |S_N|}
        \,\le\,
        {|a_{N-1}| \over |S_N|} \,
           {|t| e^{-N\gamma} \over 1 - (1+|t|)e^{-N\gamma}}$
\end{itemize}
In particular, if $N \ge [\log (1+2|t|)]/\gamma$, we have
$\Delta_N \le |a_{N-1}|$ and $\delta'_N \le |a_{N-1}|/|S_N|$.
\end{proposition} 


Let us conclude by estimating the size of the largest term
$\max\limits_n |a_n|$.
Define
\be
   b_n  \;=\; {|t|^n \, |x|^{n(n+1)/2}  \over
               (1-|x|)(1-|x|^2) \cdots (1-|x|^n)}
   \;,
\ee
so that $|a_n| \le b_n$ (with equality if $|t|=1$ and $0<x<1$).
Suppose that $|x| = e^{-\gamma}$;
it then follows from the computation in \reff{eq.lemma3.5}
that $b_n$ attains its maximum value at
$n = \lfloor \log(1+|t|) / \gamma \rfloor$,
and that this maximum value is
$\exp[C(|t|)/\gamma + O(1)]$ where
\be
   C(t)  \;=\;  \smhalf \log(1+t) \log\!\left( {t \over 1+t} \right)
                \,-\, {\rm Li}_2\!\left( {1 \over 1+t} \right)
                \,+\, {\pi^2 \over 6}
   \;.
\ee
In particular, $C(1) = \pi^2/12$ [from \reff{eq.dilog.half}].
Therefore, for $|t|=1$ the largest term can be as large in magnitude
as $e^{\pi^2/12\gamma}$ (and is indeed of this order when $0 < x < 1$);
while the answer $R(t,x)$ can be as small in magnitude as $e^{-\pi^2/6\gamma}$
(and is indeed of this order when $t=1$ and $0 < x < 1$).
It is therefore necessary to maintain,
in intermediate stages of the calculation, approximately
$(\pi^2/4\gamma)/\log 10 \approx 1.07/\gamma$ digits of
working precision beyond the number of significant digits desired
in the final answer.

\section{Comparison with other algorithms}   \label{sec4}

Let us conclude by briefly comparing the algorithm based on \reff{eq_qseries}
with some alternative algorithms for computing $R(t,x)$.

Direct use of the defining product \reff{def_R}
manifestly gives an algorithm that is only linearly convergent,
and in which the convergence rate deteriorates linearly as $|x| \uparrow 1$.
Moreover, there is severe loss of numerical precision
when multiplying numbers that are very near 1.
An alternative approach can be based on the logarithmic variant \reff{eq_logR};
this sum is again only linearly convergent,
but the problem of loss of numerical precision
is alleviated by use of the logarithm.

A slight improvement to the algorithm based on \reff{def_R}
can be obtained by noting that
\be
   \prod\limits_{n=N+1}^\infty  (1 - tx^n)
   \;=\;
   1 \,-\, {t x^{N+1} \over 1-x} \,+\, O(x^{2N})
 \;,
 \label{eq.prod.correction}
\ee
so that correcting the product \reff{def_R}
by the factor $1 - t x^{N+1}/(1-x)$
yields an estimate with error $O(x^{2N})$ rather than $O(x^N)$.
But the basic inefficiencies of the elementary algorithm remain.

Gatteschi \cite{Gatteschi_69} has proposed the following iterative algorithm
for computing $R(t,x)$:\footnote{
   I have altered his notation to conform to that of the present paper:
   his $a,q,\xi,x_n,y_n$ correspond to my $tx,x,1-\sigma,\alpha_n,\beta_n$.
   Gatteschi's algorithm has been employed by
   Allasia and Bonardo \cite{Allasia_80}.
}
Choose a complex number $\sigma \notin \{0,tx\}$
and define
\begin{subeqnarray}
   \alpha_0  & = &  1    \\[2mm]
   \beta_0   & = &  {\sigma \over \sigma - tx}  \\[4mm]
   \alpha_{n+1}  & = &  \alpha_n \, {\sigma\alpha_n + (1-\sigma)\beta_n
                                     \over \beta_n}
       \\[2mm]
   \beta_{n+1}   & = &  \alpha_n \, {\sigma\alpha_n + (1-\sigma)\beta_n
                                     \over   x \alpha_n + (1-x) \beta_n}
 \label{def.gatteschi}
\end{subeqnarray}
Gatteschi proves that
$\lim_{n\to\infty} \alpha_n = \lim_{n\to\infty} \beta_n = R(t,x)$.
In fact, it can easily be shown by induction that
\begin{subeqnarray}
   \alpha_n  & = &  \prod\limits_{k=1}^n (1-tx^k)  \\[2mm]
   \beta_n   & = &  {\sigma  \over \sigma - tx^{n+1}} \, \alpha_n
 \label{eq.solution.gatteschi}
\end{subeqnarray}
(though Gatteschi does not note this);
so the iteration \reff{def.gatteschi} gives
simply a disguised way of computing the defining product \reff{def_R}
and a slight variant of it.
Now, it is easily seen that
\begin{subeqnarray}
   {\alpha_n \over R(t,x)}  & = &  1 \,+\, {t x^{n+1} \over 1-x} \,+\, O(x^{2n})
 \slabel{eq.expansion.gatteschi.a}
      \\[2mm]
   {\beta_n \over R(t,x)}  & = &
     1 \,+\, t x^{n+1} \left( {1 \over 1-x} + {1 \over \sigma} \right)
       \,+\, O(x^{2n})
\end{subeqnarray} 
Therefore, if we set $\lambda = 1 + \sigma/(1-x)$,
the linear combination
\be
   \widehat{\alpha}_n \;\equiv\; \lambda \alpha_n + (1-\lambda)\beta_n
   \;=\;
   \left[ 1 \,-\, {\sigma tx^{N+1} \over (1-x)(\sigma - tx^{N+1})} \right]
   \alpha_n
\ee
converges to $R(t,x)$ more rapidly than either $\alpha_n$ or $\beta_n$ does
(as Gatteschi observes in a special case):
namely, $\widehat{\alpha}_n / R(t,x) = 1 + O(x^{2n})$.
But this is essentially equivalent (modulo higher-order terms)
to the ``improved'' elementary algorithm based on the correction factor
\reff{eq.prod.correction}.

Finally, Slater \cite{Slater_54,Slater_66} has computed $R(t,x)$ using
the ``other'' Euler formula\footnote{
   For a proof of \reff{eq_qseries_2},
   see e.g.\ \cite[p.~19, Corollary 2.2]{Andrews_98}
   or \cite[p.~34, Lemma 4(b)]{Knopp_70}.
}
\be
   {1 \over R(t,x)}  \;\equiv\;  \prod\limits_{n=1}^\infty (1 - tx^n)^{-1}
     \;=\;  \sum_{m=0}^\infty
            {t^m x^m \over (1-x)(1-x^2) \cdots (1-x^m)}
   \;.
 \label{eq_qseries_2}
\ee
But this algorithm is only linearly convergent;
it is no better than the logarithmic sum \reff{eq_logR},
and indeed is somewhat inferior due to the potentially small denominator.

\section*{Acknowledgments}

I wish to thank George Andrews and Mireille Bousquet-M\'elou
for suggesting (independently) that I use Euler's formula
\reff{eq_qseries} to compute $R(t,x)$;
George Andrews and Henry McKean for useful comments on
$q$-series, $q$-products and modular forms;
and Jes\'us Salas for constant close collaboration in testing
the numerical algorithm.
I also wish to thank Thomas Prellberg for drawing my attention
to his paper \cite{Prellberg_95} and to that of Moak \cite{Moak_84}.

\end{document}